\documentclass[a4paper,twoside,10pt,fleqn]{article}

\usepackage{graphics}
\usepackage{latexsym}
\usepackage{amsmath}
\usepackage{amssymb}
\usepackage{epsf}

\newcommand{\commentaar}[1]{{}}

\newcommand{\cN}{\mathcal{N}}
\newcommand{\cO}{\mathcal{O}}
\newcommand{\eps}{\varepsilon}
\newcommand{\fR}{\mathbb{R}}
\newcommand{\glob}{\textrm{glob}}

\newcommand{\loc}{\textrm{loc}}
\newcommand{\vol}{\textrm{vol}}
\newcommand{\Hess}{\textrm{Hess}}

\newtheorem{theorem}{Theorem}
\newtheorem{proposition}[theorem]{Proposition}
\newtheorem{lemma}[theorem]{Lemma}

\newcommand{\versie}{29/nov/2008}
\newcommand{\kortetitel}{Convergence}
\newcommand{\korteauteur}{Hoveijn}

\title{On convergence of the optimization process in Radiotherapy treatment planning}
\author{I. Hoveijn\\University of Groningen, Department of Mathematics\\
P.O. Box 800, 9700 AV Groningen, The Netherlands}
\date{\versie}

\begin{document}
\maketitle

\begin{abstract}\noindent
The Radiotherapy treatment planning optimization process based on a quasi-Newton algorithm with an object function containing dose-volume constraints is not guaranteed to converge when the dose value in the dose-volume constraint is a critical value of the dose distribution. This is caused by finite differentiability of the dose-volume histogram at such values. A closer look near such values reveals that convergence is most likely not at stake, but it might be slowed down.
\end{abstract}

\section{Introduction}\label{sec:intro}
A recent development in Radiotherapy treatment planning is to use an optimization process to obtain a treatment plan. To do so one needs to formalize the clinical objectives in a so called object function whose optimum should correspond to an optimal treatment plan. Many current treatment planning systems rely on a quasi-Newton algorithm to perform the optimization. In order to guarantee that this process indeed converges to an optimum the object function should at least be twice continuously differentiable. This last property is the main subject of this study. As it turns out the object functions used in Radiotherapy do not always have this property so that convergence of the optimization process can not directly be guaranteed in those cases.

In section \ref{sec:optbg} we present some background for the use of optimization in Radiotherapy treatment planning. The object functions currently used are built from so called dose-volume and EUD constraints. By dose-volume constraints, the dose-volume histogram enters the object function. It is the finite differentiability of the former that causes the insufficient differentiability of the latter. In section \ref{sec:objdiff} this relationship is explained in more detail. A example illustrates how finite differentiability of the object function influences the optimization process. The practical consequences thereof are the subject of section \ref{sec:conc}. Where section \ref{sec:objdiff} is kept informal, section \ref{sec:proof} provides formal statements and proofs of the differentiability results. This last section does not contain new arguments but serves as a support for the exposition.

\section{Background for optimization in Radiotherapy}\label{sec:optbg}
%
\subsection{Optimization in Radiotherapy}\label{sec:optrt}
In Radiotherapy treatment planning a careful choice of values for treatment parameters is necessary to obtain a dose distribution in the patient that meets criteria for eradicating a malignant tumour. At the same time it should meet criteria so that healthy tissues and organs are not critically damaged. Due to the fact that dose is deposited in every irradiated part of the patient, these criteria are conflicting in almost every clinical situation. Experienced planners are usually very good at intuitively balancing such conflicting goals, still there is a need for more objective and quantifiable methods to devise and evaluate treatment plans. This has led to the construction of object functions, depending on treatment parameters, which formalize the objectives of dose on tumours and healthy tissues, see \cite{bra,tbn,wm}. By optimizing this object function with respect to treatment parameters, one hopes to get an optimal treatment plan. Another, very different, reason to use optimization in treatment planning is the emergence of IMRT where the sheer number of treatment parameters simply precludes manual planning\commentaar{refs}. 

Even if we know that a function has a unique optimum, there is no simple recipe to find it. Therefore we have to resort to numerical, iterative methods to solve optimization problems, which, given a starting point, generate a sequence converging to an optimum. The latter can be guaranteed only if certain conditions on the object function are met. It seems, however, that if these conditions are considered at all, it is taken for granted that they are fulfilled.

Many treatment planning systems use a quasi-Newton method, like the LBFGS algorithm, to find an optimum of the object function, see for example \cite{blnz} or \cite{sb} as a general reference. One of the conditions guaranteeing that this algorithm yields a sequence locally converging to an optimum is that the object function is twice continuously differentiable. The object functions currently used in treatment planning systems are limited to so called \emph{EUD constraints} and \emph{dose-volume constraints}, or combinations of both, see for example \cite{tbn}, but also see section \ref{sec:objrt} remark \ref{rem:compro}. Our results show that the former do satisfy the differentiability condition. The latter however, do not, that is not for all values of the treatment parameters.

\subsection{Object functions in Radiotherapy}\label{sec:objrt}
The goal of Radiotherapy treatment planning is to devise a treatment plan achieving a high \emph{tumour control probability} (TCP) combined with a low \emph{normal tissue complication probability} (NTCP). A priori both TCP and NTCP are unrelated to treatment parameters. Even relating them to criteria on the dose distribution in the patient is a non-trivial task. To simplify matters one does not consider the 3-dimensional dose distribution, but a derived quantity namely the dose-volume histogram for each relevant region. The latter are tumours and healthy organs for which the dose-volume histogram tells which part of the volume receives at least a certain dose, see equation (\ref{eq:volfun}) below. Data from previously treated patients are used to relate TCP and NTCP to dose-volume histograms, see \cite{ehkms}. At this point further simplifications are necessary, because TCP and NTCP can only be related to certain points of dose-volume histograms or even only to averages. Generally for tumours this leads to dose-volume constraints, that is a minimum and a maximum dose for volumes $v_{\min}$ and $v_{\max}$ respectively, or in a more concise notation $(d_{\min},v_{\min})$ and $(d_{\max},v_{\max})$. In case of healthy tissues the situation is much more complicated, in practice one uses a mixture of dose-volume constraints and EUD constraints, see \cite{tbn}.

Criteria for the dose distribution or the derived dose-volume histograms are related to treatment parameters, thus the treatment planning objectives can be expressed with an object function of these treatment parameters. Now the object function used in treatment planning is a weighted sum of 'local' object functions, one for each dose-volume constraint or EUD constraint. Let us consider a tumour with a minimum and a maximum dose-volume constraint. Furthermore $V_{\sigma}$ is the dose-volume histogram of the tumour for a given dose distribution and treatment parameters $\sigma$, that is $V_{\sigma}(h)$ is the (relative) volume of the tumour, region $R$, receiving at least dose $h$
\begin{equation}\label{eq:volfun}
V_{\sigma}(h) = \vol \{x \in R\;|\; f_{\sigma}(x) \geq h\}.
\end{equation}
In this expression $\vol(A)$ is the volume of the region $A$ and $f_{\sigma}$ is the dose distribution. One constructs what we will call \emph{local object functions}. For example the local object function for a minimum does-volume constraint for a specific volume is a function $F_{\min}(\sigma)$ which increases the more the actual relative volume at $d_{\min}$ is below $v_{\min}$ and decreases or even becomes zero if the actual relative volume at $d_{\min}$ is above $v_{\min}$. Similarly a maximum dose object function $F_{\max}(\sigma)$ is constructed such that it increases when the relative volume at $d_{\max}$ is above $v_{\max}$ and decreases or becomes zero otherwise. Thus the object function for the tumour alone would be the weighted sum of two local object functions, usually of the following form
\begin{equation}\label{eq:objfun}
F(\sigma) = w_1\, F_{\min}(\sigma) + w_2\, F_{\max}(\sigma) = w_1\, G(v_{\min}-V_{\sigma}(d_{\min})) + w_2\, G(V_{\sigma}(d_{\max})-v_{\max}),
\end{equation}
where $w_1>0$ and $w_2>0$ are weights and
\begin{equation}\label{eq:locobjfun}
G(x) = \left\{\begin{array}{ll}
0,   & x < 0\\
x^2, & x \geq 0,
\end{array}\right.
\end{equation}
Although the object function $F$ is a function of the treatment parameters $\sigma$ it does not depend on them directly, but only via the dose-volume histograms. The latter depend on the dose distribution which directly depends on the treatment parameters.

The EUD constraints are in a similar way expressed in the \emph{equivalent uniform dose}. A quantity used to measure the biological effect of radiation especially for healthy tissues, see \cite{nie}\commentaar{nog andere refs?}. We denote the equivalent uniform dose by $E_{\alpha}$, here it is a function of treatment parameters $\sigma$
\begin{equation}\label{eq:eud}
E_{\alpha}(\sigma) = \bigg[\frac{1}{\vol(R)} \int_R f_{\sigma}(x)^{\alpha} \, dx\bigg]^{1/\alpha}.
\end{equation}

\textbf{Remarks}
\begin{enumerate}\topsep 0pt\itemsep 0pt
\item \label{rem:compro} The word \emph{constraint} in this context is somewhat misleading because a dose-volume nor an EUD constraint is used as a constraint in the sense of constrained optimization. If they were, many if not most optimization problems in treatment planning would have no solution. It is a common clinical situation that a tumour is located adjacent to a critical organ. For several reasons the region, called PTV, on which the treatment dose is prescribed is larger than the tumour. See \cite{icru} for a systematic description of XTV's, where X = G, C, I or P. Therefore the critical organ is possibly partly located in the PTV. Suppose for example that one third of the volume of the critical organ is located inside the PTV. Furthermore suppose that the minimal dose-volume constraint on the PTV is 60 Gy (Gy is unit of dose) and the maximum average dose constraint on the organ is 20 Gy. Since the dose distribution varies continuously in space and its gradient is bounded, this constrained optimization problem will have no solution. Due to the anatomy a minimum dose of 60 Gy on the PTV implies an average dose well above 20 Gy on the organ and vice versa. A weighted sum of local object functions is a way to find a clinically acceptable compromise for this situation. To avoid the term \emph{constraint} some authors use the term \emph{soft constraint}, but in my opinion this term is superfluous because we already have the notion of object function. Nevertheless we adhere to the word constraint, in order to stay in line with the existing literature.\commentaar{zij het tandenknarsend}
\item \label{rem:eud} The equivalent uniform dose, as defined in equation (\ref{eq:eud}), can be regarded as a generalized average, in fact for $\alpha=1$ it is the standard average of the dose distribution over the region $R$. Formally it is identical to the $L_p$-norm of the dose distribution $f_{\sigma}$ restricted to a region $R$ of the patient. But $p$ is called $\alpha$ and it has a biological interpretation. One associates a value of $\alpha$ to each tissue type, but there is no agreement yet in the literature about these values\commentaar{refs}.\commentaar{ Since the dose distribution is positive there is no need to restrict to $\alpha > 0$. Then there are two limiting cases, $\lim_{\alpha \to \infty} E_{\alpha}(\sigma)=\max_{x \in R}f_{\sigma}(x)$ and $\lim_{\alpha \to -\infty} E_{\alpha}(\sigma)=\min_{x \in R}f_{\sigma}(x)$. Both numbers exist because $f_{\sigma}$ is a smooth function and $R$ is compact.}\commentaar{moet dit er in? nee}
\item \label{rem:dvh} The dose-volume histogram refers to the graph of a function that is usually given the same name. This function, defined for each $\sigma$ in equation (\ref{eq:volfun}), can be regarded as a volume function in the sense of \cite{hvn}. Let us fix $\sigma$ for the moment. The dose-volume histogram of a region $R$ is a function $h \mapsto V_{\sigma}(h)$ whose value is the volume of $R$ receiving at least dose $h$. In other words it gives the volume of a region enclosed by two level sets of the dose distribution $f_{\sigma}$, namely those $x \in R$ where $f_{\sigma}(x) = h$ and the $x\in R$ where $f_{\sigma}(x) = \max_R(f_{\sigma})$ (this last level set will be a point).\commentaar{The results of \cite{hvn} imply that this function is smooth as a function of $h$ when $h$ is a regular value of $f_{\sigma}$ but only finitely differentiable when $h$ is a critical value of $f_{\sigma}$.}\commentaar{crypto!}
\item \label{rem:treatpar} Any quantity that influences the dose distribution may be considered as a \emph{treatment parameter}. Usually treatment parameters are \emph{beam intensities} and \emph{beam angles}, that is parameters of the treatment unit.
\item \label{rem:planpar} The object function $F(\sigma)$ in equation (\ref{eq:objfun}) also depends on $d_{\min}$, $d_{\max}$, $v_{\min}$, $v_{\max}$, $w_1$ and $w_2$. Such parameters might be called \emph{planning parameters}. Here we will keep them fixed.
%
%
\end{enumerate}

\section{Differentiability of object functions in Radiotherapy}\label{sec:objdiff}
The numerical optimization process based on a quasi-Newton method converges locally to an optimum if the object function is at least twice continuously differentiable. In the sequel we will simply say differentiable. In current Radiotherapy practice the object function is a function of treatment parameters, but only indirectly via dose-volume constraints and or EUD constraints. Here we will consider their differentiability properties.

Let us start with object functions only containing EUD constraints. Then $F(\sigma)$ is a weighted sum of local object functions involving only
\begin{displaymath}
G(d_{\min}-E_{\alpha}(\sigma))\;\; \text{and}\;\; G(E_{\alpha}(\sigma)-d_{\max})
\end{displaymath}
like in equations (\ref{eq:objfun}) and (\ref{eq:locobjfun}). In this case the object function $F$ is differentiable. Here we have to assume that the dose distribution is differentiable.

However, when the object function contains dose-volume constraints it is not always differentiable. This is due to the finite differentiability of the dose-volume histogram at critical values of the dose distribution, even when the latter is differentiable. Since it is not immediately obvious how the differentiability of the dose-volume histogram as a function of dose affects the differentiability of the object function as a function of treatment parameters, we present the following explanation. Essentially, the treatment parameters enter the argument of the volume function because in the presence of parameters, the critical values are parameter dependent.

When treatment parameters vary, the dose distribution $f_{\sigma}$ varies and in particular its critical values. Thus a critical value may pass through $h$. Suppose $x_{\sigma}=0$ is a critical point for $\sigma=\sigma_0$ (a critical point can always be translated to $0$). It can be shown, see section \ref{sec:proof}, that for values of $\sigma$ near $\sigma_0$ and $x$ near $0$, $f_{\sigma}$ can be transformed to a new function $\hat{f}_{\sigma}(x) = g(x) + f_{\sigma}(x_{\sigma})$. Where $g$ is a $\sigma$ independent local standard of $f_{\sigma}$ with critical point at zero and critical value zero and $f_{\sigma}(x_{\sigma})$ is the critical value of $\hat{f}_{\sigma}$ for $\sigma$ near $\sigma_0$.

Using the previous we will show in section \ref{sec:proof} that $V_{\sigma}(h) = V_{\glob}(\sigma,h) + V_{\loc}(\sigma,h)$ can be split into a global and a local part. The global part depends differentiably on $\sigma$. Let us consider the local part in more detail. On a neighbourhood $\cO \subset \fR^3$ of $0$ we have
\begin{align*}
V_{\loc}(\sigma,h) &\propto \vol \{x \in \cO \;|\; \hat{f}_{\sigma}(x) \geq h\}\\
                   &= \vol \{x \in \cO \;|\; g(x) + f_{\sigma}(x_{\sigma}) \geq h\}\\
                   &= \vol \{x \in \cO \;|\; g(x) \geq h - f_{\sigma}(x_{\sigma})\}\\
                   &= V_g(h - f_{\sigma}(x_{\sigma})).
\end{align*}
The function $V_{\loc}(\sigma,h)$ is proportional to the volume in the expression above and the proportionality factor is a differentiable function of $\sigma$. The last line in the formula above shows the volume function $V_g$ for a region enclosed by level sets of $g$. According to \cite{hvn} this function is not twice continuously differentiable at critical values of $g$. Thus $V_{\loc}(\sigma,h)$ is not twice continuously differentiable at values of $\sigma$ for which
\begin{equation}\label{eq:hyper}
h = f_{\sigma}(x_{\sigma}).
\end{equation}
This, in turn, means that the object function $F$ is not twice continuously differentiable at values of $\sigma$ such that $h$ is a critical value of $f_{\sigma}$.

\textbf{Remarks}
\begin{enumerate}\topsep 0pt\itemsep 0pt
\item \label{rem:gnotsmo} The function $G$ in equation (\ref{eq:locobjfun}) is not differentiable at $0$. This alone would already imply that the object function in (\ref{eq:objfun}) is not differentiable. However we will assume that we are in a situation like in remark \ref{rem:compro} where we are looking for a compromise. This means we stay away from points where $G$ is not differentiable.
\item \label{rem:codim} Note that equation (\ref{eq:hyper}) defines a hyper surface $\Lambda$ in parameter space $\fR^m$. It will be $m-1$-dimensional in most of its points, therefore we need to vary only one parameter to pass through $\Lambda$.
\item \label{rem:dosesmo} Throughout we assume that the dose distribution is a differentiable function. If it were less differentiable, the volume function would also be less differentiable for both critical and regular values. 
\end{enumerate}

\subsection{Examples}\label{sec:xmpl}
The aim of this section is to provide two examples. One shows a differentiable object function, the other shows an object function that fails to be twice continuously differentiable at a certain value of the treatment parameters, here we take $\sigma=0$. In the first case the Newton iteration converges to an optimum, wheras in the second case the Newton iteration does not converge if one of the iterates becomes $0$. To avoid lengthy calculations obscuring the points we wish to illustrate we will make several simplifications.

\begin{figure}[htbp]
\setlength{\unitlength}{1mm}
\begin{picture}(100,25)(0,0)
\put(48, 0){\epsffile{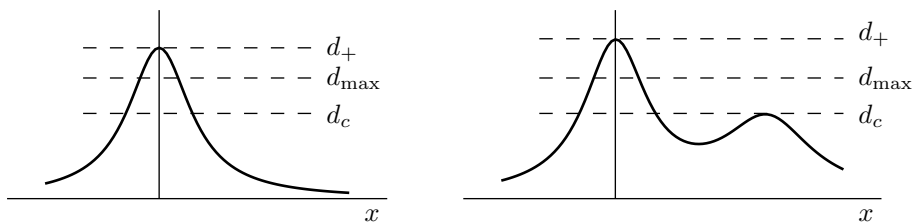}}
\put(70, 19){$d_+$}
\put(140,21){$d_+$}
\put(70, 10){$d_c$}
\put(140,10){$d_c$}
\put(70, 15){$d_{\max}$}
\put(140,15){$d_{\max}$}
\put(75, -3){$x$}
\put(140,-3){$x$}
\end{picture}
\caption{On the left the graph of $f_{\sigma_1}(x,0,0)$ and on the right the graph of $f_{\sigma_2}(x,0,0)$. $d_+$, $d_{\max}$ and $d_c$ are explained in the text.\label{fig:distroos}}
\end{figure}

Let $f_{\sigma}$ be a dose distribution on the region $R$. The object function consists of a single maximum dose-volume constraint $(d_{\max},v_{\max})$ so $F(\sigma) = G(V_{\sigma}(d_{\max}),v_{\max})$. The difference between the to examples will be the dose distribution. Suppose for $\sigma$ in a sufficiently large neighbourhood of $\sigma_1$ the dose distribution is a function with a single maximum and for a similar neighbourhood of $\sigma_2$ it is a function with two maxima for example:
\begin{displaymath}
f_{\sigma_1}(x,y,z) = \frac{1}{1+x^2+y^2+z^2}\;\;\text{and}\;\;
f_{\sigma_2}(x,y,z) = \frac{1}{1+x^2+y^2+z^2} + \frac{1}{2+(x-4)^2+y^2+z^2}.
\end{displaymath}
Here $\sigma_1$ and $\sigma_2$ are different parameter values. In figure \ref{fig:distroos} we draw the functions $f_{\sigma_1}(x,0,0)$ and $f_{\sigma_2}(x,0,0)$. For both dose distributions we call $d_+$ the value at the global maximum and $d_c$ is the value of the other local maximum of $f_{\sigma_2}$. Both $d_+$ and $d_c$ depend on the treatment parameters $\sigma$. Let us assume that always $d_c < d_{\max} < d_+$. Furthermore we assume that $V_{\sigma}(d_{\max}) < V_{\sigma}(d_c)$.

For $\sigma$ near $\sigma_1$ the dose distribution has only one critical value namely $d_+$ so for $d < d_+$, $V_{\sigma}(d)$ is a differentiable function of $\sigma$ which we call $U$. In the second case, where $\sigma$ is near $\sigma_2$, we will be particularly interested in $d=d_c$. For $d$ near $d_c$ but $d \neq d_c$, $V_{\sigma}(d)$ is again a differentiable function. At $d=d_c$, $V_{\sigma}(d) = U(\sigma) + V_g(d_c-d)$ can be split in a local and a global part, as in the previous section, where the global part is differentiable. We assume for simplicity that the global part is again $U$. The local part $V_g$ is essentially a function of the following form
\begin{displaymath}
V_g(h) = \left\{\begin{array}{ll}
\alpha h^{\frac{3}{2}},& h > 0\\
0                          ,& h \leq 0,
\end{array}\right.
\end{displaymath}
where $g$ is a local standard form of the dose distribution, $\alpha$ is a constant. Since the volume function $V_g$ is not differentiable at $0$, the object function (that we will construct shortlye hereafter) is not differentiable at values of $\sigma$ for which $d=d_c$. The latter defines a hyper surface $\Lambda$ in parameter space of codimension one, which means that we need to vary only one parameter to cross $\Lambda$. Therefore we choose new parameters and the first one we take to be $d-d_c$. Only this parameter will be relevant for the example, we call it again $\sigma$, but now $\sigma \in \fR$. Thus for the second dose distribution we finally have the following volume function:
\begin{equation}\label{eq:volfun2}\commentaar{sigma in twee betekenissen?!}
V_{\sigma}(d) = \left\{\begin{array}{ll}
U(\sigma) + \alpha (-\sigma)^{\frac{3}{2}},& \sigma < 0\\
U(\sigma)                                 ,& \sigma \geq 0.
\end{array}\right.
\end{equation}

Now we construct the object function. First we scale such that $v_{\max}=1$. Then the object function for the first dose distribution is
\begin{displaymath}
F_1(\sigma) = (U(\sigma)-1)^2.
\end{displaymath}
The object function for the second dose distribution is given by almost the same formula, but $U$ is replaced by the volume function in equation (\ref{eq:volfun2}). The function $U$ must satisfy several conditions so that it is consistent with the dose distributions we chose. First $U$ must be differentiable and decreasing, furthermore $U(0) > 1$ since we assumed that $V_{\sigma}(d_{\max}) < V_{\sigma}(d_c)$. A possible choice is $U(\sigma) = \frac{15}{10+\sigma}$.

\begin{figure}[htbp]
\setlength{\unitlength}{1mm}
\begin{picture}(100,35)(0,0)
\put(0,  15){\epsffile{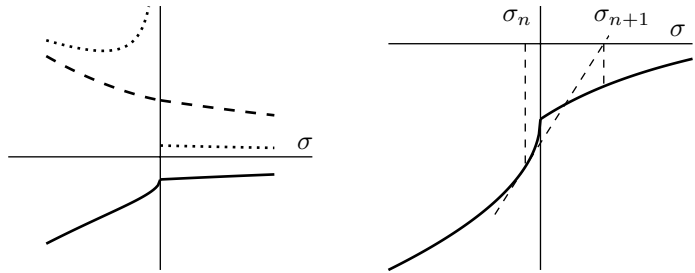}}
\put(68, 16){$\sigma$}
\put(117,31){$\sigma$}
\put(95, 33){$\sigma_n$}
\put(107,33){$\sigma_{n+1}$}
\end{picture}
\caption{In the left figure, the graph of $F_2$ is indicated by a dashed line, $F'_2$ by a solid line and $F''_2$ by dotted line. The figure on the right shows a magnification of the graph of $F'_2$ near $\sigma=0$. The dashed lines indicate how the Newton iterate $\sigma_{n+1}$ is constructed from $\sigma_n$.\label{fig:graphs}}
\end{figure}

In the first optimization problem the object function $F_1$ is differentiable so there are no difficulties in applying the Newton algortihm. However, in the second optimization problem the object function $F_2$ is not twice continuously differentiable at $\sigma=0$. Recall that the 1-dimensional Newton iteration for $F$ is $\sigma_{n+1} = \phi(\sigma_n) = \sigma_n - \frac{F'(\sigma_n)}{F''(\sigma_n)}$. Fixed points $\sigma_0$ of $\phi$ correspond to extremal points of $F$ if $F''$ is continuous and nonzero at $\sigma_0$. Since the limit from the left of the second derivative of $F_2$ is unbounded, the Newton iteration has a fixed point at $\sigma=0$ which is not related to an minimum or maximum of $F_2$. But this 'spurious' fixed point is not stable for the Newton iteration $\phi$. Only if $\sigma_n=0$ the next iterates will also be zero. Therefore it does not lead to convergence to a non extremal point of $F$ in general. Although it may slow down the convergence of the Newton iterates. For a graphical explanation see figure \ref{fig:graphs}.

\section{Conclusion}\label{sec:conc}
Our results show that the Radiotherapy treatment planning optimization process based on a quasi-Newton method locally converges to an optimum of the object function for most values of the treatment parameters. Only when the object function contains a dose-volume constraint $(h,v)$, convergence can not be guaranteed. The reason is that the object function is not sufficiently differentiable for values of the treatment parameters such that $h$ is a critical value of the dose distribution $f_{\sigma}$. Although for all other values of $\sigma$, the object function is differentiable. The question remains whether this leads to practical consequences for optimization in treatment planning.

Let us take a closer look at the values of the treatment parameters for which the object function is not sufficiently differentiable. In parameter space these values of $\sigma$ are determined by the equation $h = f_{\sigma}(x_{\sigma})$, see section \ref{sec:objdiff}, expressing that $h$ is a critical value of the dose distribution. Geometrically this equation defines a hyper surface $\Lambda$ in parameter space. Roughly speaking $\Lambda$ is $m-1$-dimensional so it has measure zero. This means that given a starting point for optimization, the Newton iteration scheme generates a sequence of points that will lie on $\Lambda$ with probability zero. The points on $\Lambda$ are spurious fixed points of the Newton iteration like in the example of the previous section. The reason is that we can generalize this example. First all types of extreme points of the dose distribution, namely minima, saddles and maxima lead to the same type of non-differentiability of the volume function. Second, also in higher dimensions the second derivative of the object function becomes unbounded at points of insufficient differentiability. Third, the spurious fixed points of the Newton iteration scheme are unstable. Our final conclusion is that in practice the convergence of the optimization process is most likely only slowed down near $\Lambda$.

\section{Proof of differentiability statements}\label{sec:proof}
In order to make a statement about the differentiability of the object function, we need a few definitions. Let $\Sigma \subset \fR^m$ be an open set of $m$ \emph{treatment parameters}. Let $D$ be a domain in $\fR^3$ representing the patient and let $R \subset D$ represent a tumour region or a healthy organ. For each $\sigma \in \Sigma$, the \emph{dose distribution} is represented by a function $f_{\sigma} : D \to \fR$ (here we do not consider $f_{\sigma}$ as a distribution in the sense of generalized functions). A function is called \emph{smooth} if it is infinitely many times continuously differentiable. We now assume that $f_{\sigma}$ is a smooth function of position $x \in D$ and moreover smoothly dependent on the parameters $\sigma \in \Sigma$. \emph{Critical points} of $f_{\sigma}$ are points in $D$ where the gradient of $f_{\sigma}$ vanishes, other points are called \emph{regular points}. The value $f_{\sigma}$ takes at a critical point is called a \emph{critical value}. If the level set $\{x \in R\;|\; f_{\sigma}(x) = h\}$ contains no critical points, $h$ is called a \emph{regular value}. We now state our results on the differentiability of object functions containing EUD constraints and dose-volume constraints respectively.

\begin{theorem}\label{the:diff1}
Let the object function $F : \Sigma \to \fR$ be given by $F(\sigma) = G(E_{\alpha}(\sigma),h)$, where $G : \fR \times \fR \to \fR$ is a smooth function and $E_{\alpha}(\sigma)$ is the equivalent uniform dose. Then $F$ is smooth with respect to $\sigma$.
\end{theorem}

\begin{theorem}\label{the:diff2} Let the object function $F : \Sigma \to \fR$ be given by $F(\sigma) = G(V_{\sigma}(h),v)$, where $G : \fR \times \fR \to \fR$ is a smooth function and $V_{\sigma}$ is a volume function for dose distribution $f_{\sigma}$. Then $F$ is smooth with respect to $\sigma$ provided that $h$ is not a critical value of $f_{\sigma}$. If $h$ is a critical value of $f_{\sigma_0}$, then $F$ is not twice continuously differentiable at $\sigma = \sigma_0$.
\end{theorem}

The proof of theorem \ref{the:diff1} is rather straightforward, but the proof of theorem \ref{the:diff2} is more involved. Since in the setting of theorem \ref{the:diff2} the differentiability of the object function $F$ is completely determined by the volume function, the proof is about the latter only. The theorem contains two statements which we here state and prove separately. The proofs heavily rely on the results in \cite{hvn}, but due to the presence of parameters we have to make several adjustments. Let us start with some definitions.

\textbf{Definitions}
\begin{enumerate}\topsep 0pt\itemsep 0pt
\item Let $f_{\sigma}: \fR^n \to \fR$ be a positive smooth function for all $\sigma \in \fR^m$ and moreover $f_{\sigma}$ depends smoothly on $\sigma$. For each $\sigma \in \Sigma$, $f_{\sigma}$ vanishes at infinity in the following sense, for each $\eps >0$ there is a compact $K \subset \fR^n$ such that for all $x \in \fR^n \setminus K$, $f_{\sigma}(x) < \eps$. We assume that critical points of $f_{\sigma}$ are non-degenerate, that is if $x$ is a critical point of $f_{\sigma}$ then $\det(\Hess\;f_{\sigma}(x)) \neq 0$.
\item Also $f: \fR^n \times \fR^m \to \fR$ is a smooth function with non-degenerate critical points. Although we consider $f$ and $f_{\sigma}$ as different functions, their values are taken to be identical: $f(x,\sigma)=f_{\sigma}(x)$ for all $x \in \fR^n$ and $\sigma \in \fR^m$.
\item We denote the level sets of $f_{\sigma}$ by $N_{h,\sigma} = \{x \in \fR^n \;|\; f_{\sigma}(x)=h\}$.
\item The level sets of $f$ are denoted by $\cN_{h} = \{(x,\sigma) \in \fR^n \times \fR^m \;|\; f(x,\sigma)=h\}$.
\end{enumerate}

\textbf{Remark.} Note that critical points of $f$ are also critical points of $f_{\sigma}$, but the opposite is not necessarily true.

\textbf{Proof of theorem \ref{the:diff1}.}
The function $f_{\sigma}$ is a positive and smooth function, smoothly depending on parameters $\sigma$. Then $[f_{\sigma}]^{\alpha}$ for $\alpha > 0$ is also smooth in both variables an parameters. Now it immediately follows that $E_{\alpha}(\sigma)$ as defined in equation (\ref{eq:eud}) is again a smooth function of $\sigma$.\hfill $\square$

The following two propositions essentially cover theorem \ref{the:diff2}.

\begin{proposition}\label{pro:diff}
Let $f_{\sigma}$ and $f$ be as defined above. Then $V_{\sigma}(h) = \vol_n \{x \in \fR^n \;|\; f_{\sigma}(x) \geq h\}$ depends smoothly on $\sigma$ at $\sigma = \sigma_0$ if $h$ is a regular value of $f_{\sigma_0}$.
\end{proposition}

\begin{proposition}\label{pro:finidiff}
Let $f_{\sigma}$ and $V_{\sigma}(h)$ be as in the previous proposition. Then $V_{\sigma}(h)$ is finitely differentiable with respect to $\sigma$ at $\sigma = \sigma_0$ if $h$ is a critical value of of $f_{\sigma_0}$.
\end{proposition}

Let us begin with the nice situation where $h$ is a regular value of $f_{\sigma_0}$, so that $V_{\sigma}(h)$ smoothly depends on $\sigma$. In the proof of proposition \ref{pro:diff} we use a property of the level sets $N_{h,\sigma}$ of $f_{\sigma}$.

\begin{lemma}\label{lem:ndiffeo}
Let $f_{\sigma}$ and $f$ be as defined above and furthermore $h$ be a regular value of $f_{\sigma_0}$. Then an open neighbourhood $\Sigma$ of $\sigma_0$ exists, such that $N_{h,\sigma}$ is diffeomorphic to $N_{h,\sigma_0}$ for all $\sigma \in \Sigma$ and the diffeomorphism depends smoothly on $\sigma$.
\end{lemma}

\textbf{Proof of proposition \ref{pro:diff}.} The volume $V_{\sigma}(h)$ is given by an integral over a region bounded by $N_{h,\sigma}$. Since, according to lemma \ref{lem:ndiffeo}, the latter depends smoothly on $\sigma$ at $\sigma=\sigma_0$, the integral and thus $V_{\sigma}(h)$ depends smoothly on $\sigma$ at $\sigma=\sigma_0$.\hfill $\square$

It remains to prove lemma \ref{lem:ndiffeo}.

\textbf{Proof of lemma \ref{lem:ndiffeo}.} Let $\{\sigma=\sigma_0\}$ be a shorthand for $\{(x,\sigma) \in \fR^n \times \fR^m \;|\; \sigma=\sigma_0\}$. Note that $N_{h,\sigma_0} = \{\sigma=\sigma_0\} \cap \cN_{h}$. Since $h$ is a regular value, $N_{h,\sigma_0}$ contains no critical points of $f_{\sigma_0}$. Thus the intersection of $\{\sigma=\sigma_0\}$ and $\cN_{h}$ is transversal. Then a neighbourhood $\Sigma$ of $\sigma_0$ exists such that for all $\sigma_1 \in \Sigma$, the intersection of $\{\sigma=\sigma_1\}$ and $\cN_{h}$ is transversal. Thus for all $\sigma_1 \in \Sigma$, $N_{h,\sigma_1}$ is a regular level set, consequently there are no critical points of $f$ in $\cN_{h,\Sigma} = \{(x,\sigma) \in \fR^n \times \fR^m \;|\; x \in N_{h,\sigma}, \sigma \in \Sigma\}$. Suppose $\sigma_1 \in \Sigma$ and $\sigma_1 \neq \sigma_0$ then an index $i$ exists such that $\sigma_{1,i} \neq \sigma_{0,i}$. Now the gradient flow of the height function $\sigma_i$ on $\cN_{h}$ is a regular flow on $\cN_{h,\Sigma}$ thus defining a diffeomorphism $\Phi$ from $N_{h,\sigma_0}$ to $N_{h,\sigma_1}$. \hphantom{knudde}\hfill $\square$

Let us now turn to the case where $h$ is a critical value of $f_{\sigma_0}$. Then there is at least one critical point $x_{\sigma_0}$ in $N_{h,\sigma_0}$, but for the sake of simplicity we assume that $x_{\sigma_0}$ is unique. In the proof of proposition \ref{pro:finidiff} we will use a similar construction as in \cite{hvn} without providing all details.

\textbf{Proof of proposition \ref{pro:finidiff}.} Let us assume that $x_{\sigma_0}$ is the unique critical point of $f_{\sigma_0}$ on the critical level set $N_{h,\sigma_0}$. Since by assumption $x_{\sigma_0}$ is a non-degenerate critical point we can apply the implicit function theorem to conclude that an open neighbourhood $\Sigma$ of $\sigma_0$ exists such that the map $\Sigma \to \fR^n : \sigma \mapsto x_{\sigma}$ and $x_{\sigma}$ is a critical point of  $f_{\sigma}$ is smooth. Now we define a new function $g_{\sigma}(x) = f_{\sigma}(x_{\sigma}+x)$, which is a smooth function of $x$ and $\sigma \in \Sigma$. Moreover $g_{\sigma}$ has a critical point at $x=0$ with critical value $f_{\sigma}(x_{\sigma})$. Then we have $V_{\sigma}(h) = \vol_n \{x \in \fR^n \;|\; g_{\sigma}(x) \geq h\}$. Following the construction in \cite{hvn} we split the latter region in a part not containing the point $0$ and a remaining small part including $0$. The volume of the first part depends smoothly on $\sigma$ using the same arguments proving proposition \ref{pro:diff}. Since $V_{\sigma}(h)$ is the sum of both volumes, the differentiability of $V_{\sigma}(h)$ is determined by that of the volume of the second part.

Let us consider this second part contained in a small neighbourhood $\cO$ of $0$. We may assume that $\cO$ is small enough to put $g_{\sigma}$ into the standard form of \cite{hvn}. In order to keep notation simple we denote the standard form again by $g_{\sigma}$, but now we can write $g_{\sigma}(x) = g(x) + f_{\sigma}(x_{\sigma})$. The actual standard form is $g(x) = \sum_{i=1}^p x_i^2 - \sum_{i=1}^q x_{p+i}^2$, with $p+q=n$. Note that $g$ no longer depends on $\sigma$, it has a critical point at $x=0$ and the critical value is $0$. This form can be obtained by a smooth change of coordinates, even smoothly depending on $\sigma$. Here again we need the non-degeneracy of critical point of $f_{\sigma}$. In the following we use the results of \cite{hvn}. As a function of $k$, the volume function $V_{\loc}(k) = \vol_n \{x \in \cO \subset \fR^n \;|\; g(x) \geq k\}$ is only finitely differentiable with respect to $k$ at $k=0$. 

Now the volume of the second part is proportional to $V_{\loc}(h-f_{\sigma}(x_{\sigma})) = \vol \{x \in \cO \subset \fR^n \;|\; g(x) \geq h-f_{\sigma}(x_{\sigma})\}$. The proportionality factor is of no concern because it originates from smooth coordinate changes and depends smoothly on $\sigma$. Thus we obtain that the differentiability of $V_{\sigma}(h)$ with respect to $\sigma$ is determined by that of $V_{\loc}(h-f_{\sigma}(x_{\sigma}))$. The latter is only finitely differentiable with respect to $\sigma$ at those values of $\sigma$ where $h-f_{\sigma}(x_{\sigma})=0$. Or, put differently, where $h$ is a critical value of $f_{\sigma}$.\hfill $\square$

\textbf{Remarks}
\begin{enumerate}\topsep 0pt\itemsep 0pt
\item The statements and proofs in this section are valid for any dimension. Therefore we have not specialized them to the case $n=3$. In this respect finitely differentiable means less than two times continuously differentiable.
\item The critical points of a parameter family of functions will generally be degenerate for certain parameter values. However we excluded this possibility because the volume function will be even less differentiable at a critical value when the critical point on the level set of this value is degenerate.
\end{enumerate}

%



\begin{thebibliography}{99}
%
\bibitem{bra} A. Brahme, Treatment Optimization using physical and radiobiological objective functions, in: "Radiation therapy physics", Alfred R. Smith (ed.) pp. 209-246, Springer Berlin, 1995.
%
\bibitem{blnz} R.H. Byrd, P. Lu, J. Nocedal, C. Zhu, A limited memory algorithm for bound constrained optimization, SIAM J. Sci. Comput. \textbf{16(5)} (1995) pp. 1190-1208.
%
\bibitem{ehkms} A. Eisbruch, R.K. Ten Haken, H.M. Kim, L.H. Marsh, J.A. Ship, Dose, volume and function relationships in parotid salivary glands following conformal and intensity-modulated irradiation of head and neck cancer, International Journal of Radiation Oncology Biology and Physics \textbf{45(3)} (1999) pp. 577-587.
%
\bibitem{hvn} I. Hoveijn, Differentiability of the volume of a region enclosed by level sets, Journal of Mathematical Analysis and Applications, \textbf{348} (2008) pp. 530-539.
%
\bibitem{icru} ICRU, Prescribing, recording and reporting photon beam therapy (supplement to ICRU report 50), ICRU Report 62, International commission on radiation units and measurements, 1999.
%
\bibitem{nie} A. Niemierko, Reporting and analyzing dose distributions: A concept of equivalent uniform dose, Medical Physics \textbf{24(1)} (1997) pp. 103-110.
%
\bibitem{sb} J. Stoer, R. Bulirsch, Introduction to  numerical analysis, Springer 1983.
%
\bibitem{tbn} C. Thielke, T. Bortfeld, A. Niemierko, From physical dose constraints to equivalent uniform dose constraints in inverse radiotherapy planning, Med. Phys. \textbf{30(9)} (2003) pp. 2332-2339.
%
\bibitem{wm} Q. Wu, R. Mohan, Algorithms and functionality of an intensity modulated radiotherapy optimization system, Med. Phys. \textbf{27(4)} (2000) pp. 701-711.
%
\end{thebibliography}
\end{document}